\documentstyle[12pt]{amsart}

\newcommand{\ext}{\mathop{ext}\nolimits}
\newcommand{\Cn}{{\Bbb C}^n}
\newcommand{\D}{\Delta}
\newcommand{\Un}{\mathop{\text{U}}\nolimits}
\newcommand{\KK}{K^{1/2}}
\newcommand{\db}{\bar\partial}
\newcommand{\Proj}{{\bf P}}
\newcommand{\R}{{\Bbb R}}
\newcommand{\im}{\sqrt{-1}}
\newcommand{\Clif}{\mathop{\text{CLIF}}\nolimits}
\newcommand{\Spin}{\mathop{\text{SPIN}}\nolimits}
\newcommand{\Or}{\mathop{\text{O}}\nolimits}
\newcommand{\SO}{\mathop{\text{SO}}\nolimits}
\newcommand{\GL}{\mathop{\text{GL}}\nolimits}
\newcommand{\E}{\mathop{\text{E}}\nolimits}
\newcommand{\gl}{\mathop{\frak{gl}\,}\nolimits}
\newcommand{\Pin}{\mathop{\rm PIN}\nolimits}
\newcommand{\spin}{\mathop{\frak{spin}\,}\nolimits}
\newcommand{\un}{\mathop{{\frak u}\,}\nolimits}
\newcommand{\so}{\mathop{\frak{so}}\nolimits}
\newcommand{\C}{{\Bbb C}}
\newcommand{\Z}{{\Bbb Z}}

\newcommand{\AAA}{{\cal A}}
\newcommand{\FF}{{\cal F}}
\newcommand{\partdel}[2]{\frac{\partial#1}{\partial#2}}
\newcommand{\liedela}[1]{L_{\partdel{}{\xi^#1}}}
\newcommand{\liedelb}[1]{L_{\overline{\partdel{}{\xi^#1}}}}
\newtheorem{thm}{Theorem}[section]

\newtheorem{prop}[thm]{Propisition}
\newtheorem{lem}[thm]{Lemma}

\theoremstyle{definition}
\newtheorem{defn}[thm]{Definition}

\numberwithin{equation}{section}

\begin{document}

\title{The Twistor correspondence of the Dolbeault complex over $\Cn$}
%complex affine space}
\author{Yoshinari Inoue}
\address{Department of Mathematics\\
	Faculty of Science\\ Kyoto University\\
	Kyoto 606-01, Japan}
\email{inoue@@kusm.kyoto-u.ac.jp}

\maketitle

\begin{abstract}
With respect to the Dolbeault complex over the flat manifold $\C^n$,
an explicit description of the inverse correspondence of the twistor
correspondence is given.

\end{abstract}

\bigbreak

\noindent
{\bf Introduction}
\bigbreak

	In [O.R], it is defined the twistor spaces of a Hermitian manifold.
Let $X$ be an $n$-dimensional Hermitian manifold and $k=0,\ldots,n$.
The $k$-th twistor space $Z_k(X)$ of $X$ is defined as a fiber
bundle over $X$ with fiber $G_{k,n}$.
The almost complex structure of $Z_k(X)$ defined in [O.R] seems
somewhat technical.
In this paper, another definition of the twistor space $Z_k(X)$ is given.
It is defined as an almost complex submanifold of the Riemannian twistor space
$Z(X)$ by using the irreducible decomposition
of the spin module $\D$ as a $U'(n)$-module. ( $\Un'(n)$ is the
double covering group of $\Un(n)$.)
$$
	\D = \bigoplus_{k=0}^n \D^k.
 $$

	In the case of Riemannian manifolds, the twistor correspondence
of the spin complex ( Dirac operator ) is given in [M].
In the Hermitian case, it can be also
defined a twistor correspondence
of the twisted Dolbeault complex:
$$
	0 \rightarrow \KK\otimes\Lambda^{0,0}
	  \buildrel \db \over \rightarrow \KK\otimes\Lambda^{0,1}
	  \buildrel \db \over \rightarrow \cdots
	  \buildrel \db \over \rightarrow \KK\otimes\Lambda^{0,n}
	  \rightarrow 0.
 $$
Let $H$ be the hyperplane bundle over $Z_k(X)$ defined by
pulling back the hyperplane bundle over $Z(X)$.
We assume the integrability of the almost complex structure of $Z_k(X)$.
An element of the cohomology group
$H^{k(n-k)}(Z_k(X), {\cal O}(H^{-n-1}))$
determines a $k$-th harmonic
form of the twisted Dolbeault complex. (Theorem \ref{z2x})

	The main theme of this paper is to give the inverse correspondence
explicitly in the case $X=\Cn$. (Definition \ref{def_x2z} and
Theorem \ref{thm_x2z})
That is, to determine an element of $H^{k(n-k)}(Z_k(X), {\cal O}(H^{-n-1}))$
from a $k$-th harmonic form of the (twisted) Dolbeault complex.
It is expected to extend this result to general Hermitian manifolds.

	Let us explain briefly the contents of this paper.
	In \S 1, we define twistor spaces of a Hermitian manifold by using
the definition of the twistor space of a Riemannian manifold in [I].
In \S 2, we describe the twistor spaces of $\Cn$ explicitly.
In \S 3, we define the inverse twistor correspondence and
prove the main theorem.

%	I would like to thank Professor K. Ueno for their encouragement
%and valuable suggestions.
%
%----------------------------------------------------------------------
%\input 1.tex
\section{The twistor spaces of a Hermitian manifold
}

By a Hermitian manifold $X$, we mean a Riemannian manifold with
a compatible almost complex structure.
In this section we define twistor spaces
 $Z_k(X), k=0,\ldots,n$ of $X$ as submanifolds of the
Riemannian twistor space $Z(X)$ defined as a
submanifold of the projectivized spinor bundle $\Proj(\D(X))$ ([I]).

	Let $\D$ be a spin module:
$$
	\D = \langle \theta_I \mid I < (1,\ldots,n) \rangle_{\bf C}
 $$
where $I<(1,\ldots,n)$ means that $I$ is a subsequence of $(1,\ldots,n)$.
For convenience, we regard  multi-indices $I,J,\ldots$ as finite sequences
of possibly duplicate elements of $\{1,\ldots,n\}$, and denote by
$IJ$ the composition of sequences $I$ and $J$.
Let us define relations among ${\theta_I}$'s as:
\begin{align*}
	  \theta_{iiI} &= - \theta_I,\\
	  \theta_{ijI} &= - \theta_{jiI}, \quad\hbox{for $i\ne j$.}
\end{align*}
Then, for any multi-index $I$,  there is a unique subsequence
$I_0$ of $(1,\ldots,n)$ such that:
$$
	\theta_I = \theta_{I_0} \quad\hbox{or}\quad \theta_I = - \theta_{I_0}.
 $$
Let $|I|$ denote the length of the reduced sequence $I_0$, and
 $i\in I$ means that the number $i$ exists in the sequence $I_0$.

With this notation, we define an action of
$\R^{2n} = \langle e_i \mid i=1,\ldots,2n \rangle_\R$ to the
spin module $\D$ as follows. For $i=1,\ldots,n$,
\begin{align*}
	  e_i \theta_I &= \theta_{iI}\\
	  e_{n+i} \theta_I &=
		\begin{cases}
		  \im \theta_{iI},	&\text{if $i \notin I$}\\
                  -\im\theta_{iI}, &\text{if $i \in I$}
		\end{cases}
\end{align*}
By the definition of the Clifford algebra,
this action can be extended to a $\Clif(\R^{2n})$-action to $\D$.
Since $\spin(2n)$ is a subspace of $\Clif(\R^{2n})$, we have
a $\spin(2n)$-action on $\D$. (See [I] for details.)

	We identify $\R^{2n}$ with $\Cn$ by the complex structure
defined by the matrix:
$$
	J = \begin{pmatrix}
		0 & -I_n\\ I_n & 0
	     \end{pmatrix}.
 $$
Then we have an inclusion $\un(n) \subset \so(2n) = \spin(2n)$.
The irreducible decomposition of $\D$ as a $\un(n)$-module
is given by:
\begin{align*}
	  \D &= \bigoplus_{k=0}^n \D^k,\\
	  \D^k &= \langle \theta_I \mid |I| = k \rangle_\C.
\end{align*}
Let $\Cn$ denote the vector representation of $\un(n)$, then we have:
\begin{equation}\label{rep_delta_k}
	\D^k \simeq (\Lambda^k \Cn) \otimes (\Lambda^n \Cn)^{-1/2}.
\end{equation}
The group $\Un'(n)$ corresponding to this representation is given by the
 following diagram:
\begin{equation*}
	\begin{array}{ccccccccc}
	0& \rightarrow& \Z_2& \rightarrow& \Pin(2n)& \rightarrow& \Or(2n)&
		\rightarrow& 0\\
	 &	      & \parallel&	 & \cup	   &		& \cup\\
	0& \rightarrow& \Z_2& \rightarrow& \Un'(n)& \rightarrow& \Un(n)&
		\rightarrow& 0.\\
	\end{array}
\end{equation*}
Let $Z$ denote the $\Pin(2n)$-orbit of $[\theta_\emptyset]\in \Proj(\D)$,
which is a complex submanifold of $\Proj(\D)$.
We define the submanifolds $Z_k, k=0,\ldots,n$ of $Z$ by
\begin{equation*}
	Z_k = Z \cap \Proj(\D^k).
\end{equation*}
This coincides with the $\Un'(n)$-orbit of
 $[\theta_{1\cdots k}]\in \Proj(\D^k)$,
hence by \eqref{rep_delta_k}, $Z_k$ is considered to be a Grassmannian
manifold $G_{k,n}$. By the definition of $Z_k$, we have a diagram:
\begin{equation}\label{zk_in_z}
	\begin{array}{ccc}
	  Z& \rightarrow&	\Proj(\D)\\
	  \uparrow&	& \uparrow\\
	  Z_k& \rightarrow& \Proj(\D^k)
	\end{array}
\end{equation}
The upper (lower) horizontal array is a $\Pin(2n)$- ($\Un'(n)$-)
equivariant mapping.
We define the hyperplane bundle $H$ over $Z_k$ by pulling back
the hyperplane bundle over $Z$ (or equivalently $\Proj(\D^k)$) .

	Let $X$ be a Hermitian manifold with spin structure:
the structure group of $X$ is $\Un'(n)$. Let $P$ denote the principal
bundle of $X$. Then we define the $k$-th twistor space and
the hyperplane bundle on it as:
\begin{align*}
	  Z_k(X) &= P \times_{\Un'(n)} Z_k\\
	      H	 &= P \times_{\Un'(n)} H.
\end{align*}
By \eqref{zk_in_z}, $Z_k(X)$ is a subbundle of $Z(X)$, and since
$Z_k$ is a complex submanifold of $Z$, it is also an almost complex
submanifold. This almost complex structure on $Z_k(X)$ is integrable
if the metric of $X$ is Bochner flat ([O.R]). In this case,
the hyperplane bundle can be considered as a holomorphic line bundle.

	Now we can define the twistor correspondence.
Let $X$ be a Hermitian manifold with spin structure and assume
that the almost complex structure of $Z_k(X)$ is integrable.
By Serre duality and \eqref{rep_delta_k}, we define:
\begin{align}\label{tw_cor}
	  T: H^{k(n-k)}(Z_k(X), {\cal O}(H^{-n-1})) \rightarrow&
	  \Gamma\left(X,\bigcup_{x\in X}
		H^{k(n-k)}(Z_k(X)_x,{\cal O}(H^{-n-1}))\right)\\
	  =& \Gamma(X, \KK\otimes \Lambda^{0,k}).\nonumber
\end{align}
\begin{thm}\label{z2x}
	An image of $T$ is a $k$-th harmonic form of the twisted
Dolbeault complex on $X$:
\begin{equation}\label{dol_cpx}
	0 \rightarrow \KK\otimes\Lambda^{0,0}
	  \buildrel \db \over \rightarrow \KK\otimes\Lambda^{0,1}
	  \buildrel \db \over \rightarrow \cdots
	  \buildrel \db \over \rightarrow \KK\otimes\Lambda^{0,n}
	  \rightarrow 0.
\end{equation}
\end{thm}

\begin{pf}
It is proved in the same way as in [H] and [M].
\end{pf}

	The main result of this paper is to give
an explicit description of the inverse of the above correspondence in
the case $X= \Cn$.
%It is expected that this result is
%extended to a general Hermitian manifold.

%\begin{conj}
%	Any $k$-th harmonic form of the twisted Dolbeault complex $(1-11)$
%is in the image of $T$.
%\end{conj}

%----------------------------------------------------------------------
%\input {header}

\section{Twistor spaces of $\Cn$}

	The Riemannian twistor space of $\R^{2n}$ is given in [I] \S 8.
Using that result, we give an explicit description of the
twistor spaces of $\Cn$.

	Let $\D'$ be a spin module of $\Spin(2n+2)$:
\begin{equation}
	\D' = \langle \theta_I \mid I< (0,1,\ldots,n) \rangle_\C.
\end{equation}
The twistor space of the $2n$-dimensional sphere is identified
with the orbit $Z'$ of the $[\theta_\emptyset]\in\Proj(\D')$
under the action of $\Pin(2n+2)$.
Since  stereographic projection defines the conformal embedding
$\R^{2n} \subset S^{2n}$, $Z(\R^{2n})$ is an open submanifold of $Z(S^{2n})$.
Let $(Z^I)$ be the homogeneous coordinates with respect to $(\theta_I)$.
Then, we have:
\begin{equation*}
	Z(\R^{2n}) = \{ (Z^I)_{I<(0,\ldots,n)} \in Z' \mid
			\exists I < (1,\ldots,n) \quad\hbox{such that}
			\quad Z^I \not= 0 \}.
\end{equation*}
Since the spin bundle of $\R^{2n}$ is trivial, we have
\begin{equation}\label{trivialization}
	Z(\R^{2n}) \simeq \R^{2n}\times Z.
\end{equation}
Since a translation of $\R^{2n}$ is a conformal transformation, it defines a
holomorphic transformation of $Z(\R^{2n})$,
 which is representable by an element of $\Spin(2n+2;\C)$.
Let $t = (x^1,\ldots,x^{2n})$ be an element of $\R^{2n}$.
Then the corresponding element of $\Spin(2n+2;\C)$ is:
$$
	\alpha(t) = 1 + \frac{1}{2} (x^1e_1 + \cdots + x^ne_n +
		x^{n+1}e_{n+2}+\cdots + x^{2n}e_{2n+1})
		(\im e_0 + e_{n+1}).
 $$
As an element of $Z'$, a point on the fiber over $0\in\R^{2n}$ is written as
$\sum_{I\not\ni 0} Z^I\theta_I$.
Hence its image by the transformation $\alpha(t)$ is:
$$
	\alpha(t)\sum_{I\not\ni0}Z^I\theta_I =
		Z^I\theta_I +
		\im \sum_{k\in I} Z^{kI}\xi^k\theta_{0I} +
		\im \sum_{k\not\in I} Z^{kI}\overline{\xi^k}\theta_{0I}.
 $$
where $\xi^k = x^k + \im x^{n+k}$ is the complex coordinate of $t$.
It follows that the projection map to the second component
of \eqref{trivialization} is given by:
\begin{equation*}
	\begin{matrix}
	  p_2:& Z(\R^{2n})& \rightarrow& Z\\
	      &	(Z^I)_{I<(0,\ldots,n)}& \mapsto& (Z^I)_{I<(1,\ldots,n)}.
	\end{matrix}
\end{equation*}
Hence we have:
\begin{equation*}
	Z_k(\Cn) = \{ (Z^I)\in Z(\R^{2n}) \mid
			Z^I = 0 \quad\hbox{for all $I<(1,\ldots,n)$ such that
		 		 $|I|\not= k$} \}.
\end{equation*}
Furthermore, the coordinate transformation of \eqref{trivialization} is
\begin{equation}\label{cood_tr}
	Z^{0J} = \im\left(\sum_{j\in J} Z^{jJ}\xi^j +
			\sum_{j \notin J} Z^{jJ}\overline{\xi^j}\right),
	\quad \text{ for $0\notin J$.}
\end{equation}
Hence the horizontal $(1,0)$-forms of $Z_k(\Cn)$ is given by
the following proposition.
\begin{prop}\label{horiz}
Horizontal $(1,0)$-forms
on $Z_k(\Cn)\simeq\C^n\times Z_k$ at $(\xi)\times(Z^I)$ is spanned by
\begin{align*}
	  \sum_{j\notin J} Z^{jJ} \overline{d\xi^j},& \quad |J|=k-1,\\
	  \sum_{j\in J} Z^{jJ} d\xi^j,& \quad |J|=k+1.
\end{align*}
\end{prop}
\begin{pf}
A point $z$ of the space $Z_k$ represents a complex structure of $\R^{2n}$.
Let us compare it with the original one corresponding to
the point $z_0 = [\theta_\emptyset]\in Z$.
For convenience, we call $(1,0)_z$-forms the $(1,0)$-forms with
respect to the complex structure corresponding to $z$.
Put $z_1=\theta_{(12\cdots k)}$.
Then $(1,0)_{z_1}$-forms are spanned by
$\overline{d\xi^1},\ldots,\overline{d\xi^k},d\xi^{k+1},\ldots,d\xi^n$.
This mean that the complexified form $\omega$ is of type $(1,0)_{z_1}$
if and only if the both $(1,0)_{z_0}$-part and $(0,1)_{z_0}$-part of
$\omega$ are of type $(1,0)_{z_1}$.
Hence the proposition follows  immediately from \eqref{cood_tr}.
\end{pf}

%---------------------------------------------------------------------------
%\input header

\section{Main Theorem}

	By restricting the $\SO(2n)$-action on $Z$, we have
a $\Un(n)$-action on $Z_k$. This action can be complexified and
we have a $\GL(n;\C)$-action on $Z_k$. This action defines a
$\C$-linear mapping:
\begin{equation*}
	\FF: \gl(n;\C) \rightarrow \Gamma(Z_k,\Theta),
\end{equation*}
where $\Theta$ is the holomorphic tangent bundle on $Z_k$,
which is considered to be the $(1,0)$-part of the complexified
tangent bundle $TZ_k\otimes\C$.
With respect to Lie algebra structures, we have:
\begin{equation}\label{lie_alg}
	\FF([a,b]) = - [\FF(a),\FF(b)].
\end{equation}
Let $({E^j}_i)$ denote a standard basis of $\gl(n;\C)$.
We define a vector field ${\FF^j}_i$ to be:
\begin{equation*}
	{\FF^j}_i = - \FF({}^t{\E^j}_i).
\end{equation*}
We define operators acting on the differential forms on $Z_k(\C^n)$:
\begin{align*}
	  D_\alpha &= \sum \ext(d\xi^a) i(\overline{{\FF^a}_b})
			L_{\partdel{}{\xi^b}},\\
	  D_{\beta}&= \sum L_{\overline{\partdel{}{\xi^a}}}
			i(\overline{{\FF^a}_b}) \ext(\overline{d\xi^b}),
\end{align*}
where $\ext$ and $i$ denote the exterior and inner multiplication,
 and $L$ denotes the Lie derivative.

Let $F$ be a power series defined as:
$$
	F(x) = \sum_{k=0}^{\infty} \frac{x^k}{(k!)^2}
 $$
%Then its derivatives are:
%$$
%	F^{(l)}(x) = \sum_{k=0}^{\infty} \frac{x^k}{k! (k+l)!}
% $$
This function and its derivatives play an important role by the
following property.
\begin{lem}\label{diff_eq}
Let $l$ be a non-negative integer. We have an equality
\begin{equation*}
	xF^{(l+2)}(x) + (l+1)F^{(l+1)}(x) - F^{(l)}(x) = 0.
\end{equation*}
\end{lem}

	Let $\Lambda_V^{0,k(n-k)}$ denote the line subbundle of
$\Lambda^{0,k(n-k)} Z_k(\Cn)$ spanned by vertical forms.
If we identify $H^{-1}$ with $\overline H$ by the standard Hermitian metric,
we have
$$
	\Lambda^{0,k(n-k)}\otimes H^{-n-1} \supset
		\Lambda_V^{0,k(n-k)}\otimes H^{-n-1}
		\simeq \overline H
 $$
since we have $\Lambda_V^{0,k(n-k)} \simeq {\overline H}^{\,-n}$.
On the other hand, we have an isomorphism
$\Gamma(Z_k,{\cal O}(H)) \simeq K^{-1/2}\otimes\Lambda^{k,0}$.
Hence we have a natural map:
$$
	j: \Gamma(\Cn, \KK\otimes\Lambda^{0,k}) \rightarrow
		\Gamma(Z_k(\Cn), \Lambda_V^{0,k(n-k)}\otimes H^{-n-1})
 $$

\begin{defn}\label{def_x2z}
We define an operator:
$$
\begin{matrix}
	\AAA:& \Gamma(\Cn, \KK\otimes\Lambda^{0,k}) &\rightarrow
		&\Gamma(Z_k(\Cn), \Lambda^{0,k(n-k)}\otimes H^{-n-1})\\
		& a &\mapsto
		& k! (n-k)! F^{(k)}(D_\beta)F^{(n-k)}(D_\alpha)j(a).
\end{matrix}
 $$
\end{defn}

Now we can state the main theorem.

\begin{thm}\label{thm_x2z}
Let $f$
% = \sum f_{\overline I}\overline{d\xi^I}$
be a twisted $(0,k)$-form on $\C^n$.
Then $\AAA(f)$ is $\db$-closed if and only if $f$ is a harmonic form
with respect to the twisted Dolbeault complex \eqref{dol_cpx}.
The restriction of $\AAA$ to the space of harmonic forms is the
inverse correspondence of $T$ of~\eqref{tw_cor}
\end{thm}

%Let us define few more operators as follows:
With respect to the standard complex structure of $\C^n$, we have
a decomposition of the complexified horizontal cotangent bundle:
\begin{equation*}
	\Lambda^1_H\otimes\C = \Lambda^{(1,0)}_H \oplus \Lambda^{(0,1)}_H,
\end{equation*}
By this decomposition, we define two projections as follows:
$$
	\pi_\gamma : \Lambda^1 \rightarrow
	  \begin{cases}
	    \Lambda^{(1,0)}_H, & \gamma = \alpha,\\
	    \Lambda^{(0,1)}_H, & \gamma = \beta.
	  \end{cases}
 $$
Let $\gamma$ be $\alpha$ or $\beta$. We define:
\begin{equation*}
	d_\gamma = \pi_\gamma\circ d.
\end{equation*}

\begin{lem}\label{oprts}
\begin{enumerate}
\item
Let $\gamma$ be $\alpha$ or $\beta$, and put $E_\gamma = [d,D_\gamma]$.
Then, we have
\begin{align*}
	E_\alpha &= - \ext(d\xi^a) L_{\overline{{\FF^a}_b}}
			 L_{\partdel{}{\xi^b}},\\
	E_\beta &= L_{\overline{\partdel{}{\xi^a}}}
			L_{\overline{{\FF^a}_b}}
			\ext(\overline{d\xi^b}).
\end{align*}
\item
 Let $f$ be a power series and $\gamma$ be $\alpha$ or $\beta$, then
\begin{equation*}
	[d, f(D_\gamma)]
	  = f'(D_\gamma)E_\gamma - f''(D_\gamma)D_\gamma d_\gamma.
\end{equation*}
\item\label{gamma}
Put $\Gamma = [E_\beta, D_\alpha]$. Then, we have
\begin{equation*}
	\Gamma  = \ext(d\xi^a)i(\overline{{\FF^a}_b})\ext(\overline{d\xi^b})
                  \liedela{c}\liedelb{c}
                + \liedelb{a}i(\overline{{\FF^a}_b})\liedela{b}
                  \ext(d\xi^c\wedge\overline{d\xi^c}),
\end{equation*}
$$
	\Gamma D_\alpha = D_\alpha \Gamma,
 $$
$$
	[E_\beta, f(D_\alpha)] = f'(D_\alpha)\Gamma.
 $$
\end{enumerate}
\end{lem}

\begin{pf}(1) and (3) are proved by simple computation.\newline
(2) It suffices to prove
$$
	[d,D_\gamma^n] =
		 n D_\gamma^{n-1}E_\gamma - (n-1)nD_\gamma^{n-1}d_\gamma,
 $$
which can be proved by induction on $n$ by using the formula:
\begin{equation}
\label{ee}
	[E_\gamma, D_\gamma] = -2D_\gamma d_\gamma.
\end{equation}
By (1), we have:
$$
	[E_\alpha,D_\alpha] =
		- \ext(d\xi^a\wedge d\xi^c)
			i([\overline{{\FF^a}_b},\overline{{\FF^c}_d}])
		  L_{\partdel{}{\xi^b}}L_{\partdel{}{\xi^d}}
 $$
Since, by \eqref{lie_alg}, we have:
$$
	[{\FF^a}_b,{\FF^c}_d] =
		\delta^a_d {\FF^c}_b - \delta^c_b{\FF^a}_d
 $$
hence we have shown \eqref{ee} when $\gamma = \alpha$.
The case $\gamma = \beta$ can be proved in a similar way.
\end{pf}

By Lemma \ref{diff_eq} and Lemma \ref{oprts}, we have:
\begin{align*}
	dF^{(k)}(D_\beta)F^{(n-k)}(D_\alpha)
	 = &F^{(k)}(D_\beta)F^{(n-k)}(D_\alpha)(d - d_\alpha - d_\beta)\\
	 &  + F^{(k)}(D_\beta)F^{(n-k+1)}(D_\alpha)
	       \{E_\alpha + (n-k+1)d_\alpha\}\notag\\
	 &  + F^{(k+1)}(D_\beta)F^{(n-k)}(D_\alpha)
	       \{E_\beta + (k+1)d_\beta\}\notag\\
	 &  - F^{(k+1)}(D_\beta)F^{(n-k+1)}(D_\alpha)\Gamma\notag
\end{align*}
Since $j(f)$ is harmonic in the vertical direction,
$$
	F^{(k)}(D_\beta)F^{(n-k)}(D_\alpha)(d - d_\alpha - d_\beta)j(f) = 0,
 $$
and by Lemma \ref{oprts} \eqref{gamma}, if $f$ is harmonic, we have
$$
	\Gamma j(f) \equiv 0 \quad\text{ modulo $(1,0)$-forms.}
 $$

To compute the action of $E_{\alpha}$ and $E_{\beta}$,
we have to take local coordinates of $Z_k$.
Let $I$ be a subsequence of $(1,\ldots,n)$ of length $k$. Then
$$
	w_{ij}=Z^{ijI}/Z^I,\quad i\in I, j\not\in I
 $$
are local coordinates of
$U_I=\{(Z^J)_{J<(1\cdots n)}\in Z_k \mid Z^I\ne 0\}$.
Put $z^J = Z^J/Z^I$ for $J<(1,\ldots,n)$.
\begin{lem}\label{vectorfield}
The vector field ${\FF^a}_b$ is written in the local coordinates as:
\begin{equation*}
	{\FF^a}_b = - \sum_{\begin{Sb}i\in I\\ j\not\in I\end{Sb}}
			z^{aiI}z^{bjI}\partdel{}{w_{ij}}.
\end{equation*}
\end{lem}
Let $\delta$ be a function:
$$
	\delta(a) =
		\begin{cases}
		  1,& \text{if $a$ is true,}\\
		  0,& \text{otherwise.}
		\end{cases}
 $$
The following lemma is a translation of the Pl\"ucker relation.
\begin{lem}\label{dirty}
  \begin{enumerate}
     \item We fix a multi-index $I$ and local coordinates $w_{ij}$ as above.
      \begin{equation*}
	\frac{\partial z^J}{\partial w_{ij}} =
	  \begin{cases}
	    -z^{ijJ}& \text{if $i\not\in J$, $j\in J$,}\\
	    0& \text{otherwize.}
	  \end{cases}
      \end{equation*}
    \item Let $J$, $K$ be multi-indices of length $|J|=k+1$, $|K|=k-1$.
      \begin{equation*}
	\sum_{a\in J\setminus K} Z^{aJ}Z^{aK} = 0
      \end{equation*}
    \item Let $J$, $K$ be multi-indices of length $|J| = |K| = k$.
      \begin{equation*}
	\sum_{j\in J\setminus K} Z^{ajJ}Z^{bjK} =
 	  - \delta(a\not\in K) Z^JZ^{baK}
	  + \delta(a\not\in J) Z^{baJ} Z^K.
      \end{equation*}
\end{enumerate}
\end{lem}
\begin{pf}
  (1) It is obvious if $|J\setminus I| \le  1$. If $|J\setminus I| \ge 2$,
    let $i_k\in I\setminus J$, $j_k\in J\setminus I$, $k = 1,2$, be numbers
   such that
    $i_1 \ne i_2$, $j_1 \ne j_2$. Then, by [I] Corollary 3.3, we have
    \begin{equation*}
      Z^J = \frac{1}{Z^{Ji_1i_2j_1j_2}}
	(-Z^{Ji_2j_2}Z^{Ji_1j_1} + Z^{Ji_2j_1}Z^{Ji_1j_2}).
    \end{equation*}
    Hence we complete the proof by induction.\newline
  (2) This follows by induction on $|J\setminus K|$ by using (1).\newline
  (3) We prove the case in which $a\ne b$. Put
\begin{align*}
	A &= \sum_{j\in J\setminus K} Z^{ajJ}Z^{bjK}\\
	  &= \sum_{j\in J\setminus K} Z^{jaJ}Z^{jbK}
		- 2\delta(a\in J\setminus K) Z^JZ^{baK}
		- 2\delta(b\in J\setminus K) Z^{abJ}Z^K.
\end{align*}
    Since
\begin{align*}
	P = (J\setminus K)\setminus (aJ \setminus bK) &=
	  \delta(a\in J\setminus K)a \cup \delta(b\in J\setminus K)b,\\
	M = (aJ\setminus bK)\setminus (J \setminus K) &=
	  \delta(a\not\in J\cup K)a \cup \delta(b\in J\cap K)b,
\end{align*}
we have:
\begin{align*}
	A &= \sum_{j\in aJ\setminus bK} Z^{jaJ}Z^{jbK}
		+ \sum_{j\in P} Z^{jaJ}Z^{jbK}
		- \sum_{j\in M} Z^{jaJ}Z^{jbK}
		-2 \sum_{j\in P} Z^{jaJ}Z^{jbK}\\
	  &= - \sum_{j\in P\cup M} Z^{jaJ}Z^{jbK}\qquad \text{by (2)}\\
	  &= - \delta(a\not\in K) Z^J Z^{baK}
		+ \delta(a\not\in J) Z^{baJ} Z^K.
\end{align*}
The case in which $a=b$ can be proved in the same way.
\end{pf}

\begin{lem}\label{very_dirty}
Let $I$ be a subsequence of $(1,\ldots,n)$ and $\overline{s^I}$
be the image of $\Lambda_{i\in I} \overline{d\xi^i}$ by~$j$.
\begin{align}
\label{eq1}
	L_{\overline{\FF^a}_b} \overline{s^I} &=
	    -(k+1)\delta^a_b \overline{s^I}
	    -\delta(a\notin I) \overline{s^{baI}}
	    -(n+1)\sum_{J \ni a}
		\frac{z^{baJ}\overline{z^J}}{N} \overline{s^I}\\
\label{eq2}
	L_{\overline{\FF^a}_b} \overline{s^I} &=
	     (n-k+1)\delta^a_b \overline{s^I}
	    +\delta(b\in I)\overline{s^{abI}}
	    +(n+1)\sum_{J \not\ni b}
		\frac{z^{abJ}\overline{z^J}}{N} \overline{s^I}
\end{align}
where
$$
	N = \sum_J |\frac{z^J}{z^I}|^2
 $$

\end{lem}
\begin{pf}
%These formula are proved by straightforward computation.
Let $w_{ij}$ be the local coordinates as above.
Let $\overline{\rho^I}$ be the image of the standard trivialization
of $\overline{H}$
over $\C^n\times U_I$ by the isomorphism $\overline{H}\simeq H^{-1}$.
Let $\overline{K^I}$ be the standard trivialization of $\Lambda^{0,k(n-k)}_V$
over $\C^n\times U_I$. Then
\begin{equation*}
	\overline{s^I} =
	 {\overline{\rho^I}}^{\otimes(n+1)}\otimes \overline{K^I}.
\end{equation*}
First, we compute:
\begin{equation*}
\begin{align*}
	L_{\overline{{\FF^a}_b}} \overline{\rho^I} &=
		\nabla_{\overline{{\FF^a}_b}} \overline{\rho^I}\\
		&= -\frac{\overline{{\FF^a}_b}(N)}{N} \overline{\rho^I}\\
		&= \sum_{\begin{Sb}i\in I\\ j\not\in I\end{Sb}}
			\frac{\overline{z^{aiI}z^{bjI}
				\frac{\partial}{\partial w_{ij}} N}}
			     {N}
			\overline{\rho^I}
			\qquad\text{[By Lemma \ref{vectorfield}]}\\
		&= - \sum_{J}\sum_{\begin{Sb}
					i\in I\setminus J\\
					j\in J\setminus I
				   \end{Sb}}
			\frac{z^J\overline{z^{ijJ}}\overline{z^{aiI}}
				\overline{z^{bjI}}}
			     {N}
			\overline{\rho^I}
			\qquad\text{[By Lemma \ref{dirty} (1)]}\\
		&= - \sum_J \sum_{i\in I\setminus J}
			\frac{z^J\overline{z^{aiI}z^{biJ}}}{N}
			\overline{\rho^I}
			\qquad\text{[By Lemma \ref{dirty} (3)]}\\
		&= ( \sum_{J \not\ni a}
			\frac{z^J\overline{z^{baJ}}}{N}
			- \delta(a\not\in I)\overline{z^{baI}} )
			\overline{\rho^I}
			\qquad\text{[By Lemma \ref{dirty} (3)]}
\end{align*}
\end{equation*}
By changing indices we have:
\begin{equation}\label{line}
	L_{\overline{{\FF^a}_b}} \overline{\rho^I} =
	\begin{cases}
		(-\delta^a_b - \delta(a\not\in I)\overline{z^{baI}}
		 - \sum_{J\ni a} \frac{z^{baJ}\overline{z^J}}{N})
		\overline{\rho^I},& \text{for \eqref{eq1}}\\
		(\delta^a_b + \delta(b\in I)\overline{z^{abI}}
		 + \sum_{J\not\ni b} \frac{z^{abJ}\overline{z^J}}{N})
		\overline{\rho^I},& \text{for \eqref{eq2}}
	\end{cases}
\end{equation}
Second, we compute:
\begin{equation*}
	L_{\overline{{\FF^a}_b}} d\overline{w_{ij}} =
		d\overline{{\FF^a}_b(w_{ij})} =
		(-\delta^i_b\overline{z^{aiI}} +
		  \delta^a_j\overline{z^{bjI}})d\overline{w_{ij}} + \cdots
\end{equation*}
Hence:
\begin{equation*}
\begin{align*}
	L_{\overline{{\FF^a}_b}} \overline{K^I} &=
		\sum_{\begin{Sb} i\in I\\ j\not\in J \end{Sb}}
		(-\delta^i_b\overline{z^{aiI}} +
		 \delta^a_j\overline{z^{bjI}}) \overline{K^I}\\
	&= \{-(n-k)\delta(b\in I)\overline{z^{abI}} +
	     k\delta(a\not\in I)\overline{z^{baI}} \} \overline{K^I}.
\end{align*}
\end{equation*}
By changing indices we have:
\begin{equation}\label{ext}
        L_{\overline{{\FF^a}_b}} \overline{K^I} =
        \begin{cases}
		\{(n-k)\delta^a_b + n\delta(a\not\in I)\overline{z^{baI}}\}
			\overline{K^I},& \text{for \eqref{eq1}},\\
		\{-k\delta^a_b - n\delta(b\in I)\overline{z^{abI}}\}
			\overline{K^I},& \text{for \eqref{eq2}}.
        \end{cases}
\end{equation}
Hence, by \eqref{line} and \eqref{ext}, we complete the proof.
\end{pf}

The next lemma follows immediately from Proposition \ref{horiz} and
Lemma \ref{very_dirty}.
\begin{lem}
Let $f_{\overline I}\overline{s^I}$ be an image of $\AAA$. Then we have
\begin{align*}
	\{E_\beta + (k+1)d_\beta\}f_{\overline I}\overline{s^I} &\equiv
	   - \sum_{a \not\in I}
	       \partdel{f_{\overline I}}{\overline{\xi^a}}
		\overline{d\xi^b}\wedge\overline{s^{baI}},\\
	\{E_\alpha + (n-k+1)d_\alpha\}f_{\overline I}\overline{s^I} &\equiv
	   - \sum_{b \in I}
	       \partdel{f_{\overline I}}{\xi^b}
		d\xi^a\wedge \overline{s^{abI}},
\end{align*}
where $\equiv$ denotes the equivalence modulo $(1,0)$-forms.
\end{lem}

This lemma shows that the coefficient of $\overline{d\xi^b}\overline{s^I}$
( $d\xi^a\overline{s^I}$ )
in $\AAA(f)$ is equal to the coefficient of $\overline{s^{bI}}$
( $\overline{s^{aI}}$ )
in the image of the Dirac operator.
Therefore, we complete the proof of the theorem.

%\input tail
%-----------------------------------------------------------------------

\end{document}